

\documentstyle{amsppt}
\TagsOnRight
\hoffset=.1in

\define\sh{\operatorname{sh}}
\define\ch{\operatorname{ch}}

\define\adx{\operatorname{ad}}

\define\ld{\ldots}

\define\a{\alpha}
\define\be{\beta}
\define\vf{\varphi}
\define\lm{\lambda}
\define\ve{\varepsilon}
\define\gm{\gamma}
\define\de{\delta}
\define\si{\sigma}
\define\om{\omega}

\define\vd{\varDelta}
\define\vl{\varLambda}
\define\vg{\varGamma}
\define\vo{\varOmega}

\define\vt{\varTheta}

\define\cb{\Bbb C}

\define\caa{\Cal A}

\define\car{\Cal R}

\define\xti{\widetilde x}
\define\vdti{\widetilde \varDelta}

\define\pd#1#2{\dfrac{\partial#1}{\partial#2}}

\define\vc#1{(#1_1,\ldots,#1_n)}
\define\vct#1{[#1_1,\ldots,#1_n]}
\define\vect#1{\{#1_1,\ldots,#1_n\}}


\define\jak{\dfrac{i}{\kappa}}
\define\ijak{\dfrac{1}{\kappa}}
\define\jakk{\dfrac{i}{\kappa^2}}
\define\ijakk{\dfrac{1}{\kappa^2}}
\define\mak{\dfrac{P_0}{\kappa}}
\define\makk{\dfrac{P_0}{2\kappa}}
\define\maks{{P_0}\slash{\kappa}}
\define\kap{{\Cal P}_\kappa}
\define\kapti{\widetilde{\Cal P}_\kappa}
\define\vep{\vec{P\,}^2}
\define\kub{\dfrac{\vec{P\,}^2}{2\kappa}}
\define\kubik{\dfrac{\vec{P\,}^2}{2\kappa^2}}


\define\1{$\kappa$-Poincar\'e group}
\define\2{bicovariant}
\define\3{calculus}
\define\4{differential}
\define\5{dimension}
\define\6{quantum}
\define\7{coefficient}
\define\8{operator}
\define\9{properties}
\define\0{algebra}

\define\irt{ \vartriangleright\!\blacktriangleleft }
\define\tri{ \blacktriangleright\!\vartriangleleft }
\define\trir{ \triangleright }
\define\tril{ \triangleleft }

\define\kapo{$\kappa$-Poincar\'e}

\define\linv{left-invariant}
\define\rinv{right-invariant}
\define\rep{representation}
\define\bic{bicrossproduct}
\define\three{three-dimensional}

\topmatter
\title The bicovariant differential calculus\\
on the three-dimensional $\kappa$-Poincar\'e group
\endtitle
\rightheadtext{Differential calculus}
\author Piotr Kosi\'nski$^*$, Micha\l \/ Majewski$^*$\\
{\it Department of Theoretical Physics}\\
{\it University of \L \'od\'z}\\
{\it ul. Pomorska 149/153, 90--236 \L \'od\'z, Poland}\\
Pawe\l \/ Ma\'slanka$^*$\\
{\it Department of Functional Analysis}\\
{\it University of \L \'od\'z}\\
{\it ul. St. Banacha 22, 90--238 \L \'od\'z, Poland}
\endauthor
\leftheadtext{P. Kosi\'nski, M. Majewski, P. Ma\'slanka}
\thanks
*\ \ \ Supported by KBN grant  2 P 302\,217\,06\,p\,02\newline
\endthanks

\abstract The \2 \4 \3 on the three-\5al $\kappa$-Poincar\'e group and the
corresponding Lie-\0 structure are described. The  equivalence of this
Lie-\0 structure and the three-\5al $\kappa$-Poincar\'e \0 is proved.
\endabstract

\endtopmatter

\document

\head I. Introduction
\endhead

Recently, considerable interest has been paid to the deformations of groups
and \0s of space-time symmetries [1].  In particular, an interesting
deformation of the Poincar\'e \0 [2] as well as group [3] has been
introduced which depend on dimensionful deformation parameter $\kappa$; the
relevant objects are called \kapo{} \0 and \1, respectively. Their structure
was studied in some detail and many of their \9 are now well understood. The
\kapo{} \0 and group for the space-time of any \5 has been defined [4], the
realizations of the \0 in terms of \4 operators acting on commutative 
Min\-kowski as well as momentum spaces were given [5]; the unitary \rep{s} 
of the deformed group were found [6]; the deformed universal covering 
$ISL(2,\cb)$ was constructed [7]; the bicrossproduct [8] structure, both of 
the \0 and group was revealed [9]. The proof of formal duality between \1 
and \kapo{} \0 was also given, both in two [10] as well as in four \5s [11].
 One of the important problems is the construction of the \2 \4 \3 on \1. 
Using an elegant approach due to Woronowicz [12], the \4 calculi on 
four-\5al Poincar\'e group [13], as well as on the Min\-kowski space [14] 
were constructed.

In the present paper we briefly sketch the construction of the  \4 \3 on the 
three-\5al \1. Apart from possible applications to the three \5al field
theory the \3 presented here is very interesting on its own and 
significantly differs from the calculi defined on four-\5al ([13]) and 
two-\5al ([15]) \1.

It is well known that in most cases of \2 \4 calculi on the quantum groups 
the space of (say) \linv{} 1-forms has larger \5ality than its classical 
counterpart. However, in many cases it is sufficient to add one additional 
biinvariant form; the corresponding \linv{} vector field reduces, in the 
classical limit, to quadratic Casimir operator [15], [16]. In the case under 
consideration it appears that it is necessary to add further invariant form 
related to Pauli-Lubanski invariant.

The paper is organized as follows. The  remaining part of the introduction 
is devoted to the description of \kapo{} \0 and group and their formal 
duality. In Section II the \2 $*$-\3 on three-\5al \1 is constructed. In 
Section III we obtain the corresponding Lie \0 and prove its equivalence to 
the \kapo{} \0. Finally, in Section IV some conclusions are given. Some 
calculations are relegated to the Appendix.

The three-\5al \1 $\kap$ is the Hopf $*$-\0 defined as follows [3]. Consider 
the universal $*$-\0 with unity, generated by selfadjoint elements 
$\vl^\mu{}_\nu$, $x^\mu$ subject to the following relations
$$
\aligned
&[x^\mu,x^\nu]  = \dfrac{i}{\kappa}( \de_0{}^\mu x^\nu -  \de_0{}^\nu 
x^\mu),\\
&[ \vl^\mu{}_\nu, x^\rho]  = - \dfrac{i}{\kappa} ((\vl^\mu{}_0 - 
\de_0{}^\mu) \vl^\rho{}_\nu + (\vl^0{}_\nu - \de_\nu{}^0)  g^{\mu\rho})
\endaligned
\tag{1.1}
$$
here $g_{\mu\rho} = g^{\mu\rho} = $ diag\,$(+, -, -)$ is the  metric 
tensor.

The  comultiplication, antipode and counit are defined as follows
$$
\aligned
&\vd(\vl^\mu{}_\nu)  =  \vl^\mu{}_\rho \otimes \vl^\rho{}_\nu,\\
&\vd(x^\mu)  = \vl^\mu{}_\nu \otimes x^\nu + x^\mu \otimes I,\\
&S(\vl^\mu{}_\nu) =  \vl_\nu{}^\mu,\\
& S(x^\mu) = - \vl_\nu{}^\mu x^\nu,\\ 
&\ve(\vl^\mu{}_\nu)  =  \de^\mu{}_\nu,\\
&\ve(x^\mu) = 0.
\endaligned
\tag{1.2}
$$
It is easy to see [9] that $\kap$ has the \bic{} structure [8]
$$
\kap = T^* \irt C(S0(2,1))
\tag{1.3}
$$
where $C(S0(2,1))$ is the standard Hopf \0 of functions defined over the 
Lorentz group, while $T^*$ is defined by the relations
$$
\aligned
&[x^\mu,x^\nu]  = \dfrac{i}{\kappa}( \de_0{}^\mu x^\nu -  \de_0{}^\nu 
x^\mu),\\
&[ \vd(x^\mu) = x^\mu \otimes I + I \otimes x^\mu,\\
& S(x^\mu) = -x^\mu,\\
& \ve(x^\mu) = 0.
\endaligned
\tag{1.4}
$$
Indeed, it is sufficient to define the following structure functions
$$
\aligned
& \be(x^\mu) = \vl^\mu{}_\nu \otimes x^\nu,\\
& \vl^\mu{}_\nu \tril x^\rho = [\vl^\mu{}_\nu, x^\rho].
\endaligned
\tag{1.5}
$$

The \three{} \kapo{} \0 $\kapti$ was first introduced in the third paper of 
[2]. We present it below in the Majid and Ruegg basis ([9]). It is a 
quantized universsal envelopping \0 in the sense of Drinfeld ([18]) 
described by the following relations
$$
\aligned
& [P_\mu,P_\nu] = 0,\\
& [M,P_0] = 0,\\
& [M,P_k] = i \ve_{kl} P_l,\\
& [M,N_k] = i \ve_{kl} N_l,\\
& [N_1,N_2] = -iM,\\
& [N_i,P_j] = i\de_{ij} \Big( \dfrac{\kappa}{2}(1 - e^{-2\maks}) + 
\dfrac{1}{2\kappa} \vep\Big) - \dfrac{i}{\kappa} P_iP_j,\\
& \vd P_0 = P_0  \otimes I + I \otimes P_0,\\
& \vd P_i = P_i  \otimes e^{-\maks} + I \otimes P_i,\\
& \vd M_i = M_i  \otimes I + I \otimes M_i,\\
& \vd N_i = I \otimes N_i + N_i \otimes e^{-\maks} + \dfrac{1}{\kappa} 
\ve_{ij} M \otimes P_j,\\
&S(P_0) = -P_0,\\
&S(P_i) = -e^{\maks} P_i,\\
&S(M_i) = -M_i,\\
&S(N_i) = - N_i e^{\maks} + \dfrac{1}{\kappa} \ve_{ij} M \, P_j e^{\maks},\\
&\ve(X)  = 0, \qquad, X = P_\mu, M,N_i,
\endaligned
\tag{1.6}
$$
here $\mu = 0,1,2$ and $i,j,k = 1,2$.

Again, we can write
$$
\kapti = T \tri U(so(2,1))
\tag{1.7}
$$
where $U(so(2,1))$ is classical envelopping \0 of $so(2,1)$ while $T$ is 
defined as follows
$$
\aligned
& [P_\mu,P_\nu] = 0,\\
& \vd P_0 = P_0 \otimes I + I \otimes P_0,\\
& \vd P_i = P_i  \otimes e^{-\maks} + I \otimes P_i,\\
&S(P_0) = -P_0,\\
&S(P_i) = -e^{\maks} P_i,\\
&\ve(P_\mu) = 0.
\endaligned
\tag{1.8}
$$
In order to show that (1.7) holds it is sufficient to define ([9])
$$
\aligned
& M \trir P_\mu = [M,P_\mu] ,\\
& N_i \trir P_\mu = [N_i,P_\mu] ,\\
&\de(M) = M \otimes I,\\
&  \de(N_i) = N_i  \otimes e^{-\maks} + \dfrac{1}{\kappa} \ve_{ij} M \times 
P_j.
\endaligned
\tag{1.9}
$$
It has been shown recently [11] that $\kap$ and $\kapti$ are formally dual 
to each other. The relevant duality relations read
$$
\aligned
& \langle \vl^{\mu_1}{}_{\nu_1} \ld \vl^{\mu_n}{}_{\nu_n}, M_{\a\be} \rangle 
=  i \sum^n_{k=1} (\de^{\mu_k}_\a g_{\nu_k \be} - \de^{\mu_k}_\be g_{\nu_k 
\a}) \prod_{l\ne k} \de^{\mu_l}_{\nu_l},\\
& \langle : F(x^\mu) :, f(p_\nu) \rangle = f \Big( i 
\dfrac{\partial}{\partial x^\nu}\Big) F(x^\mu)_{|x=0}.
\endaligned
\tag{1.10} 
$$
The normal product $: F(x^\mu) :$ is defined as the one in which all $x^0$ 
factors stand  leftmost. The variables $x^\mu$ on the right-hand side of 
(1.10) are viewed as commuting ones.

\head II. Bicovariant $*$-calculus on \1
\endhead

Let us recall the main result of Woronowicz theory. Given a Hopf \0 $\caa$ 
and $a \in \caa$, we write
$$
(\vd \otimes I) \circ \vd(a) = (I \otimes \vd) \circ  \vd(a) = \sum_k a_k 
\otimes b_k \otimes c_k
\tag{2.1}
$$
and define the adjoint action on $a$
$$
\adx (a) =  \sum_k b_k \otimes S(a_k) c_k.
\tag{2.2}
$$
According to Woronowicz ([12]) a \2 $*$-\3 is uniquely defined by the choice 
of right ideal $\car\subset \ker \ve$ which has the following \9
\roster
\item"{(i)}" $\car$ is $\adx$-invariant, i.e. $\adx(\car) \subset \car 
\otimes \caa$,
\item"{(ii)}" for any $a \in \car$, \ $S(a)^* \in \car$.
\endroster
The standard \3 in the commutative case is obtained by choosing $\car = 
(\ker \ve)^2$. Below (Theorem 1) we construct  the relevant ideal for \1. It 
can be obtained as follows. One starts with $(\ker \ve)^2$. However, due to 
the noncommutativity, it fails to satisfy (i). It appears that this can be 
cured by adding to the generators of $(\ker \ve)^2$ some terms linear in  
generators of $\ker \ve$ (and proportional to $1\slash \kappa$); new  
generators  form a (not completely reducible) multiplet under the adjoint 
action of $\kap$. Moreover, (ii) is also satisfied. The whole procedure 
would be rather straightforward were it not for the fact that the ideal 
obtained is  
too large --- it produces the \3 which has smaller \5ality than its 
classical counterpart. Therefore, we are forced to consider a smaller ideal. 
This can be achieved by subtracting some subrepresentations from the 
representation spanned by the generators under consideration. We performed 
this subtraction in the most economical way: one trace and one completely 
antisymmetric representations have been subtracted.

The final result can be summarized as follows.

Let us introduce the following notation
$$
\aligned
&\vd^\mu{}_\nu  =  \vl^\mu{}_\nu - \de^\mu{}_\nu,\\
&\vd^\mu{}_\nu{}^\a \equiv x^\a( \vl^\mu{}_\nu - \de^\mu{}_\nu) - 
\dfrac{i}{\kappa} (\de^0{}_\nu(\vl^{\mu\a} - g^{\mu\a}) + \vl^\mu{}_0 
(\vl^\a{}_\nu - \de^\mu{}_\nu)),\\
& x^{\mu\nu} \equiv x^{\mu} x^{\nu} + \dfrac{i}{\kappa} (g^{\mu\nu} x^0 - 
g^{0\mu} x^\nu).
\endaligned
\tag{2.3} 
$$

Then the following theorem holds

\proclaim{Theorem I} Let $\car  \subset \ker \ve$ be the right ideal 
generated by the following elements
$$
\align
&\vd^\a{}_\be  \vd^\mu{}_\nu,\\
&\vdti^{\mu\nu\a} \equiv   \vd^{\mu\nu\a} - \dfrac{1}{6} \ve_{\rho\si\gm} 
\vd_{\rho\si\gm},\\
& \xti^{\mu\nu}  \equiv x^{\mu\nu} - \dfrac{1}{3} g^{\mu\nu} x^\a{}_\a.
\endalign
$$
Then $\car$ has  the following \9
\roster
\item"{(i)}" $\car$ is $\adx$-invariant,  $\adx(\car) \subset \car 
\otimes \kap$,
\item"{(ii)}" for any $a \in \car$, \ $S(a)^* \in \car$.
\item"{(iii)}" $\ker \ve\slash \car$ is spanned by the following elements
\endroster
$$
\align
&x^\mu; \ \ \vl^\mu{}_\nu, \ \ \mu< \nu; \qquad \vf \equiv x^\a{}_\a = x^2 + 
\dfrac{2i}{\kappa} x^0,\\
&  \vd \equiv  \ve_{\mu\nu\a}  \vd^{\mu\nu\a}.
\endalign
$$
\endproclaim

We shall omit the proof of this theorem which goes along the same lines as 
in the four-\5al case (see the second paper of [13]) and is long. Let us 
only note the following:
\roster
\item"{(a)}" $\vd^\mu{}_\nu{}^\a$ and $x^{\mu\nu}$ are 'improved" generators 
of $(\ker \ve)^2$ ($\vd^\a{}_\be{} \vd^\mu{}_\nu{}$ is the same as in the 
classical case) which form the (not completely reducible) multiplet under 
adjoint action of $\kap$ (see Appendix);
\item"{(b)}" the ideal generated by $\vd^\a{}_\be{} \vd^\mu{}_\nu{}$, 
$\vd^\mu{}_\nu{}^\a$, $x^{\mu\nu}$ equals $\ker \ve$; in order to obtain 
reasonable (in the sense that it contains all \4s $dx^\mu$, $d\vl^\a{}_\be$) 
\3 we have subtracted the trace of $x^{\mu\nu}$ and  completely 
antisymmetric part of $\vd^\mu{}_\nu{}^\a$;
\item"{(c)}" it is easy to conclude from (iii) that our \3 is eight-\5al.
\endroster
Having established the structure of $\car$ we can now follow closely the 
Woronowicz construction. First, we define the \linv{} 1-forms. The basis of 
the space of \linv{} 1-forms is spanned by the elements of the form 
$\pi r^{-1}(I \otimes a_i)$, where $\{a_i\}$ is a basis in $\ker \ve\slash 
\car$, \, the operation $r^{-1}$ is defined by
$$
r^{-1} (a \otimes  b) = (a \otimes  I) (S \otimes  I) \vd(b)
\tag{2.6}
$$
and the operation $\pi$ is defined by
$$
\pi (\sum a_k \otimes b_k) = \sum a_k db_k
\tag{2.7}
$$
where
$$
\kap \otimes  \kap \ni \sum a_k \otimes  b_k
$$
is such an element that
$$
\sum a_k b_k = 0.
$$
We readily obtain
$$
\aligned
&\om^\mu{}_\nu \equiv  \pi r^{-1} (I  \otimes (\vl^\mu{}_\nu - 
\de^\mu{}_\nu)) = \vl_\a{}^\mu  d  \vl^\a{}_\nu,\\
&\om^\a \equiv  \pi r^{-1} (I  \otimes x^\a) =  \vl_\mu{}^\a d x^\mu,\\
&\om \equiv \pi r^{-1} (I  \otimes \vf) = d\vf - 2x_\mu dx^\mu,\\
&\vo \equiv \pi r^{-1} (I  \otimes \vd) = \ve_{\nu\a\be} \vl_\si{}^\nu 
\om^\be \vl^{\si\a} + \dfrac{2i}{\kappa}  \ve_{o\a\be}  \om^{\a\be}.
\endaligned
\tag{2.8} 
$$
The next step is to find the commutation rules between the invariant forms 
and  the generators  of $\kap$. This is again straightforward and can be 
briefly explained as follows. 

We write
$$
b \om_i - \om_i b = \pi r^{-1} r [(b \otimes I) r^{-1} (I \otimes a_i)  - 
r^{-1} (I \otimes a_i) (I \otimes b)] 
$$
where
$$
r(a \otimes b) = (a \otimes I) \vd (b).
$$
Due to the formulae ([12])
$$
r((a \otimes I)q) = (a \otimes I)r(q), \qquad r(q(I \otimes a)) = r(q) \vd 
(a)
$$
we get
$$
r[(b \otimes I) r^{-1} (I \otimes a_i) -  r^{-1} (I \otimes a_i) (I \otimes 
b)]  = b \otimes a_i - (I\otimes a_i) \vd (b).
$$
Then we use the relations describing our ideal $\car$ to simplify the 
right-hand side and apply again $\pi \circ r^{-1}$. The detailed 
calculations results in the following formulae
$$
\aligned
& [\vl^\mu{}_\nu, \om^\a{}_\be] = 0,\\
& [x^\a, \om^\mu{}_\nu] = - \jak(\de^0_\nu \vl^\a{}_\rho \om^{\mu\rho} + 
\de^\mu_0 \vl^\a{}_\rho  \om^\rho{}_\nu - \vl^\a{}_\nu \om^\mu{}_0 - 
\vl^{\a\mu} \om^0{}_\nu)\\
& \qquad \qquad \qquad \qquad - \dfrac{1}{6} \ve^\mu{}_\nu{}^\be 
\vl^\a{}_\be \vo,\\
& [\vl^\mu{}_\nu, \om^\a] = - \jak(\de^0_\nu \vl^\mu{}_\rho \om^{\rho\a} + 
\vl^\mu{}_0   \om^\a{}_\nu)  -  \dfrac{1}{6} \ve^\rho{}_\nu{}^\a 
\vl^\mu{}_\rho \vo,\\
&[x^\mu,\om^\a] = -\dfrac{1}{3} \vl^{\mu\a} \om + \jak \vl^{\mu\a} \om^0 - 
\jak \de^\a_0 \vl^\mu{}_\rho  \om^\rho,\\
&[\vl^\mu{}_\nu,\om] = \dfrac{3}{\kappa^2} \vl^\mu{}_\rho  \om^\rho,\\
&[x^\mu,\om] =   \dfrac{3}{\kappa^2} \vl^\mu{}_\rho  \om^\rho{}_\nu,\\
&[\vl^\mu{}_\nu,\vo] = 0,\\
&[x^\mu,\vo] = \dfrac{3}{\kappa^2} \ve_{\be\rho\si} 
\vl^{\mu\be}\om^{\rho\si}.
\endaligned
\tag{2.10} 
$$
Let us now define the right action of $\kap$ on 1-forms ([12])
$$
{}_\vg \vd(adb) = \vd(a)(d \otimes I)\vd(b).
$$
Simple calculations give
$$
\aligned
&{}_\vg{\vd(\om^\mu{}_\nu)} = \om^{\rho}{}_{\si} \otimes \vl_\rho{}^\mu 
\vl^\si{}_\nu,\\
&{}_\vg{\vd(\om^\mu)} = \om^{\rho}{}_{\si} \otimes \vl_\rho{}^\mu x^\si + \om^\rho 
\otimes \vl_\rho{}^\mu ,\\
&{}_\vg {\vd(\om)} = \om \otimes I,\\
&{}_\vg{(\vo)} = \vo \otimes I
\endaligned
\tag{2.11} 
$$
which gives the following \rinv{} forms
$$
\aligned
&\eta^\mu{}_\nu = \om^\be{}_\gm \vl^\mu{}_\be \vl_\nu{}^\gm,\\
&\eta^\mu = - \om^\be{}_\gm \vl_\rho{}^\mu x^\rho \vl^\mu{}_\be + \om^\be 
\vl^\mu{}_\be,\\
&\eta = \om,\\
&\vt = \om.
\endaligned
\tag{2.12} 
$$
This concludes the description of bimodule $\vg$ of 1-forms on $\kap$. 
External \0 can be now constructed as follows [12]. On $\vg^{\otimes 2}$ we 
define a bimodule homomorphism $\si$ such that
$$
\si(\om \otimes_{{}_{\kap}} \eta) = \eta \otimes_{{}_{\kap}} \om
\tag{2.13}
$$
for any \linv{} $\om \in \vg$ and any \rinv{} $\eta \in \vg$. Then by 
definition
$$
\vg^{\wedge 2} = \dfrac{\vg^{\otimes 2}}{\ker(I - \si)}.
\tag{2.14}
$$
Higher external power of $\vg$ can be constructed in a similar way [12]. 
The result of action of $\si$ on our forms is given in Appendix. Finally, 
after long analysis we obtain the following set of relations
$$
\aligned
& \om^\mu{}_\nu \wedge \om^\a{}_\be +  \om^\a{}_\be \wedge \om^\mu{}_\nu = 
0,\\
&   \om \wedge  \om = 0,\\
& \vo \wedge  \vo = 0,\\
& \om \wedge  \vo + \vo \wedge \om = 0,\\
&  \vo \wedge \om^\mu{}_\nu + \om^\mu{}_\nu \wedge  \vo = 0,\\
&  \om^\mu{}_\nu \wedge  \om + \om \wedge \om^\mu{}_\nu - 
\dfrac{3}{\kappa^2} \om^\si{}_\nu \wedge \om_\si{}^\mu = 0,\\
& \om^\mu{}_\nu  \wedge \om^\a + \om^\a \wedge \om^\mu{}_\nu  + 
\jak \de^\nu{}_0 \om^\a{}_\rho  \wedge  \om^\mu{}_\rho + \jak \de_0{}^\mu 
\om^\a{}_\rho \wedge \om^\rho{}_\nu\\
& \qquad  -\jak \om^\a{}_\nu \wedge \om^{\mu}{}_0 - \jak \om^{\a\mu} \wedge 
\om_{0\nu} = 0,\\ 
& \om^{\mu} \wedge   \vo + \vo \wedge \om^{\mu} + \dfrac{3}{\kappa^2} 
\ve^{\mu}{}_{\rho\be} \om_\nu{}^\be \wedge \om^{\rho\nu} = 0,\\
& \om \wedge \om^\mu + \om^\mu \wedge \om - \dfrac{3}{\kappa^2} 
\om^{\mu\rho} \wedge \om_\rho = 0,\\
& \om^{\mu} \wedge \om^{\nu} + \om^{\nu} \wedge \om^{\mu} + \jak \de_0{}^\nu 
\om^{\mu\rho} \wedge \om_\rho  + \jak \de_0{}^\mu \om^{\nu\rho} \wedge 
\om_\rho = 0.
\endaligned
\tag{2.15} 
$$
The basis in $\vg^{\wedge 2}$ consists of the following elements
$$
\align
& \om^{\a\be}  \wedge \om^{\mu\nu} \ \ \ ( \a < \be, \ \    \mu < \nu, \ \   
(\a \be) \ne (\mu\nu), \ \ \a < \mu),\\
&  \om^{\a\be} \wedge \om^\mu, \ \   \om^{\a\be} \wedge \om, \ \ \om^{\a\be} 
\wedge \vo, \ \ \om^\a \wedge \om^\mu \ \ (\a < \mu),\\
& \om^\a  \wedge \om, \ \ \om^\a \wedge \vo, \ \  \om \wedge \vo
\endalign
$$
and $\dim \vg^{\wedge 2} = \binom{8}{2}$.

The next step is to introduce the $\ast$-operation. Due to Theorem I this 
can be done considently; moreover, it is sufficient to consider 1-forms. One 
gets
$$
\aligned
& (\om^\mu{}_\nu)^*  = \om^\mu{}_\nu,\\
& (\om^\mu)^* = \om^\mu - \jak \om^\mu{}_0,\\
& \om^* = -\om,\\
& \vo^* = -  \vo.
\endaligned
\tag{2.16}
$$
To complete our exterior \3 we derive the Cartan--Maurer equations. The 
existence of external derivative is granted by Theorem 4.1 of [12]. (2.8) 
imply then the following formulae for  external derivatives of \linv{} forms
$$
\aligned
& d\om^\mu{}_\nu  = \om_\rho{}^\mu \wedge \om^\rho{}_\nu,\\
& d\om^\mu = \om_\rho{}^\mu  \wedge \om^\rho,\\
& d\om = 0,\\
& d\vo = 0.
\endaligned
\tag{2.17}
$$

\head III. The Lie algebra like structure
\endhead

Let us derive the counterpart of the classical Lie \0. To this end we 
introduce the counterparts of the \linv{} fields. They are defined by the 
formula
$$
da = \dfrac{1}{2}(\chi_{{\phantom{}}_{\mu\nu}} \ast a) \om^{\mu\nu} + 
(\chi_{{\phantom{}}_{\mu}} \ast a)\om^\mu + (\chi_{{\phantom{}}_{}} \ast a) 
\om  + (\lm \ast a)  \vo
\tag{3.1}
$$ 
where, for any linear functional $\vf$ on $\kap$,
$$
\vf \ \ast \ a \equiv (I \otimes \vf) \vd(a).
\tag{3.2}
$$

The product of two functionals $\vf_1$, $\vf_2$  is defined by the standard 
duality relation
$$
\vf_1 \vf_2(a) \equiv (\vf_1 \otimes \vf_2) \vd(a).
\tag{16}
$$
In order to find the Lie \0 structure, we apply the external derivative to 
both sides of (3.1), we use  $d^2a = 0$ on the left-hand side and calculate 
the right-hand side using (2.17) and again (3.1). Nullifying the 
coefficients  in front of basis elements of $\vg^{\wedge 2}$,  we find the 
\6  Lie \0
$$
\aligned
& [m,l_i] = \Big(1 + \jak \chi_{{\phantom{}}_{0}} - \dfrac{3}{\kappa^2} 
\chi\Big) \ve_{ik}l_k - \jak \chi_{{\phantom{}}_{i}} m + \dfrac{6}{\kappa^2} 
\lm \chi_{{\phantom{}}_{i}} ,\\
& [m_i,l_j] = -\Big(1 + \dfrac{2i}{\kappa}  \chi_{{\phantom{}}_{0}} - 
\dfrac{3}{\kappa^2} \chi\Big) \ve_{ij} m + \dfrac{6}{\kappa^2} \ve_{ij} \lm
 \chi_{{\phantom{}}_{0}} ,\\
& [m, \chi_{{\phantom{}}_{0}} ] = 0,\\
& [m, \chi_{{\phantom{}}_{i}} ] = \Big(1 + \jak \chi_{{\phantom{}}_{0}} - 
\dfrac{3}{\kappa^2} \chi\Big)  \ve_{ik}  \chi_{{\phantom{}}_{k}},\\
& [l_i, \chi_{{\phantom{}}_{0}} ] = \Big(1 + \jak \chi_{{\phantom{}}_{0}} - 
\dfrac{3}{\kappa^2} \chi\Big)  \chi_{{\phantom{}}_{i}},\\
& [l_i, \chi_{{\phantom{}}_{k}} ] = \de_{i_k} \Big(1 + \jak 
\chi_{{\phantom{}}_{0}} - \dfrac{3}{\kappa^2} \chi\Big) 
\chi_{{\phantom{}}_{0}},\\ 
&[\chi_{{\phantom{}}_{\mu}}, \chi_{{\phantom{}}_{\nu}}] = 0
\endaligned
\tag{3.4}
$$
while $\lm$ and $\chi$ commute with all generators;  here $m = 
\chi_{{\phantom{}}_{12}}$, $l_i = \chi_{{\phantom{}}_{i0}}$.

Finally, the involution for $\chi$'s and $\lm$ can be defined using (2.16) 
together with the Woronowicz duality relations ([12]), which in our case, read
$$
\aligned
&\langle \om^{\mu\nu}, \chi_{{\phantom{}}_{\a\be}} \rangle = \de^\mu_\a 
\de^\nu_\be - \de^\mu_\be \de^\nu_\a,\\
&\langle \om^{\nu}, \chi_{{\phantom{}}_{\mu}} \rangle = \de^\nu_\mu ,\\
&\langle \om, \chi \rangle = 1 ,\\
&\langle \vo, \lm   \rangle = 1.
\endaligned
\tag{3.5}
$$
the reamining brackets being vanishing. If one puts
$$
\langle \psi, \eta^* \rangle = - \langle \psi^*, \eta \rangle
$$
one readily gets
$$
\aligned
&\chi^*_{{\phantom{}}_{\mu}} = - \chi_{{\phantom{}}_{\mu}},\\
& m^* = - m,\\
& l_i^* = - l_i - \jak \chi_{{\phantom{}}_{i}} ,\\
& \chi^* = \chi,\\
&\lm^* = \lm.
\endaligned
\tag{3.5}
$$
Having our \6 Lie \0 constructed, we can now pose the question what is the  
relation between our functionals and the  elements of the \kapo{} \0 
$\kapti$. 

It has been shown ([11]) that $\kap$  and $\kapti$ are formally dual. 
Therefore we expected  our functionals to be expressible in terms of the 
elements of $\kapti$.

To show this let us note that, in the notation introduced in (1.6) the 
following substitutions reproduce (3.4)
$$
\aligned
&  \chi_{{\phantom{}}_{0}} = -i \Big( \kappa \sh \big(\mak \big) + \kub 
e^{\maks} \Big),\\ 
&  \chi_{{\phantom{}}_{i}} = -i e^{\maks} P_i   ,\\
& \chi = - \dfrac{1}{6} \Big(4\kappa^2 \sh^2  \big( \makk \big)- \vep 
e^{\maks} \Big)   ,\\ 
& \lm = - \dfrac{1}{6} \Big( \Big(\kappa\sh  \big(\mak\big) + \kub e^{\maks} \Big) M 
+ (P_1N_2 - P_2N_1) e^{\maks}\Big),\\
& m =  -i e^{\maks} M - \dfrac{6i}{\kappa}  \lm ,\\
& l_i = -i e^{\maks} N_i.
\endaligned
\tag{3.7}
$$
Then, from the Woronowicz theory, it follows that  the coproducts of the 
functionals $\vf_i$ ($\vf_i = \chi_{{\phantom{}}_{\mu\nu}}, 
\chi_{{\phantom{}}_{\nu}}, \chi,\lm$) can be written in the form
$$
\vd \vf_i = \sum_j \vf_j \otimes f_{ji} + I \otimes \vf_i
\tag{3.8}
$$
where $f_{ji}$ are the functionals entering commutation rules between the 
\linv{} forms and elements of $\kap$
$$
\om_ja = \sum_i (f_{ji} \ \ast \ a) \om_i
\tag{3.9}
$$
They are therefore calculable, in principle at least, in terms of what we 
already know. In order to check (3.7) we have first to compare the 
coproducts. Using (1.6) and (3.7) one obtains
$$
\align
\vd\chi_{{\phantom{}}_{0}} & = \chi_{{\phantom{}}_{0}} \otimes 
\Big(\ch\big(\mak\big) + \kub e^{\maks} \Big) + \chi_{{\phantom{}}_{i}} 
\otimes \ijak P_i\\
& + \chi \otimes \dfrac{3i}{\kappa} \Big(\sh\big(\mak\big) + \kub e^{\maks} 
\Big) + I \otimes \chi_{{\phantom{}}_{0}} \\
\vd\chi_{{\phantom{}}_{i}} & = \chi_{{\phantom{}}_{i}} \otimes I + 
\chi_{{\phantom{}}_{0}} \otimes \ijak P_i e^{\maks} + \chi \otimes 
\dfrac{3i}{\kappa^2} P_i e^{\maks} + I \otimes \chi_{{\phantom{}}_{i}}\\
\vd\chi & = \chi_{{\phantom{}}_{0}} \otimes \Big(-\dfrac{i}{3}\Big) 
\Big(\kappa \sh\big(\mak\big) + \kub e^{\maks} \Big) + 
\chi_{{\phantom{}}_{i}} \otimes \Big(-\dfrac{i}{3}\Big) P_i\\
& + \chi \otimes \Big(\ch\big(\mak\big) + \kub e^{\maks}\Big) + I \otimes 
\chi\\
\vd\lm & = \lm \otimes \Big( \ch\big(\mak\big) - \kub e^{\maks} \Big) + 
\chi_{{\phantom{}}_{j}} \otimes \Big(-\dfrac{i}{6\kappa} P_jM - 
\dfrac{i}{6} \ve_{jk} N_k\Big) \\
& + \chi_{{\phantom{}}_{0}} \otimes \Big[- \dfrac{i}{6} \Big(\kappa\ch
\big(\mak\big) + \kub e^{\maks}\Big) M - \dfrac{i}{6\kappa} \ve_{jk} P_j 
N_k e^{\maks}\Big]\\
& + \chi \otimes \dfrac{1}{2\kappa^2}\Big[ \Big(\kappa\sh \big(\mak\big) + 
\kub e^{\maks}\Big) M  + \ve_{jk} P_j N_k e^{\maks}\Big]   \\
& l_k \otimes \Big(-\dfrac{i}{6}\Big) \ve_{jk} P_j + m \otimes 
\Big(-\dfrac{i}{6}\Big) \Big(\kappa\sh \big(\mak\big) - \kub e^{\maks}\Big)\\
& + I \otimes \lm
\endalign
$$
$$
\aligned
\vd l_i & = l_i \otimes I + \chi_{{\phantom{}}_{0}} \otimes \ijak e^{\maks} 
N_i+ \chi \otimes  \dfrac{3i}{\kappa^2}  e^{\maks} N_i + m \otimes \ijak  
\ve_{ik} e^{\maks}P_k \\  
&  + \lm \otimes \dfrac{6i}{\kappa^2} \ve_{ik} e^{\maks} P_k + I \otimes l_i\\
\vd m & = m \otimes \Big( \ch\big(\mak\big) + \kubik e^{\maks} \Big) \\
& +  \chi_{{\phantom{}}_{0}} \otimes \Big[\ijak \Big(\ch \big(\mak\big) - 
\kub  e^{\maks}\Big) M - \ijakk \ve_{jk} P_j N_k  e^{\maks} \Big] \\
& + \chi \otimes \dfrac{3i}{\kappa^3}  \Big(\Big(\kappa\ch \big(\mak\big) 
-  \kub e^{\maks}\Big) M - \ve_{jk} P_j N_k e^{\maks}\Big)\\
& + \lm \otimes \Big(\dfrac{6i}{\kappa}\Big) \Big(\sh \big(\mak\big) + 
\kubik e^{\maks}\Big) + \chi_{{\phantom{}}_{j}}  \otimes \Big(-\ijakk P_j M 
- \ijak \ve_{jk} N_k\Big)\\
& + l_k \otimes \Big(- \ijak\Big) \ve_{jk} P_j + I \otimes m.
\endaligned
\tag{3.10}
$$

Formulae (3.10) have the expected form (3.9). To show that the functionals 
$f_{ji}$ in (3.8), (3.9) coincide with those defined by (3.10), let us note 
the following. The coassociativity implies the relations ([12])
$$
\vd f_{ji} = \sum_k f_{jk} \otimes f_{ki}.
\tag{3.11}
$$
The same relations must hold true for the  functionals appearing on the 
right-hand sides of (3.10) (the coproduct for $\kapti$ is consistent and 
coassociative as well). Therefore it is suffcicient  to check that the two 
sets of functionals $f_{ji}$ defined above take the same values on 
generators of $\kap$. The relevant values for the functionals defined in 
(3.9) are readily obtained by using the explicit form of the commutation 
rules (2.10). On the other hand, for the  functionals defined by (3.10) we 
use the duality  $\kap \Longleftrightarrow \kapti$ established in [11]. As 
it was mentioned in the introduction, both  $\kap$ and $\kapti$ have the 
\bic{} structure. From the \bic{} theory it follows that the generic element 
of $\kap$ (resp. $\kapti$) can be written as $X \otimes \vl$ (resp. $P 
\otimes M$) where $X$ (resp. $P$) is an arbitrary element of $T^*$ (resp. 
$T$) while $\vl$ (resp. $M$) --- an  arbitrary element of $C(S0(2,1))$  
(resp. $U(so(2,1))$). The duality relations can be written as
$$
\langle X \otimes \vl, P \otimes M \rangle = \langle X,P \rangle \langle 
\vl,M \rangle.
\tag{3.12}
$$
Using (3.12) together with (1.10) we have checked that both definitions of 
functionals $f_{ji}$ indeed coincide. In order to complete the proof of 
(3.7) it is now sufficient to compare the values of both sides on generators 
of $\kap$. For the left-hand side we use (3.1) while the right-hand side is 
calculated with the help of  duality relations (1.10) and (3.12).

Thus we have shown that the functionals appearing in the Woronowicz 
formulation of  differential \3 are formally expressible in terms of 
elements of $\kapti$. In particular,  $\chi$ is proportional to the first 
Casimir operator, while $\lm$ provides a deformation of Pauli-Lubanski 
invariant.

As in the four-\5al case the relation between $\chi$ and 
$\chi_{{\phantom{}}_{\mu}} $ does not follow from Cartan-Maurer equations. 
However, contrary to the  four-\5al case, the relation
$$
\lm = \dfrac{1}{12} \ve^{\mu\nu\a} \chi_{{\phantom{}}_{\mu}}   
\chi_{{\phantom{}}_{\nu\a}} ,
\tag{3.13}
$$
which is readily obtainable from (3.7), is not derivable from 
Cartan-Maurer equations in contrast with its four-\5al counterpart.

\head IV. Conclusions
\endhead

We have constructed \2  $*$-calculi on three-\5al \1. The starting point 
was the Woronicz theory of \4 calculi on \6 groups. The main ingredience of 
this approach is the choice of right ideal in $\ker \ve$ which is invariant 
under the adjoint action of the group. In the  classical case the ideal 
under consideration is $(\ker \ve)^2$. In order to obtain as slight as 
possible deformation of classical \3 we have started with the generators of  
$(\ker \ve)^2$. However, it appeared that they do not form a multiplet under 
the action  of the \1. To cure this we have modified them by adding the 
appropriate, $\kappa$-dependent terms. Due to the noncommutativity, the 
ideal generated in this way coincided with the whole  $\ker \ve$. Therefore, 
it appeared necessary to subtract from the new generators some 
$\adx$-invariant terms. As a consequence the resulting \3 contains more 
invariant forms than its classical counterpart. This results in  increasing 
the \5s of the relevant  {\it Lie \0}. The additional elements provide the 
deformations of mass squared invariant and Pauli-Lubanski invariant. The 
Woronowicz theory provides us with a  unique generalization of the notion of 
\linv{} vector fields. We have shown  that the relevant functionals are 
expressible in terms of generators of \kapo{} \0, the number of the latter 
being equal to the \5 of the classical Poincar\'e \0. Therefore, there must 
be relation between the Woronowicz functionals. In our case the functionals 
$\chi$ and $\lm$ are expressible in terms of the functionals 
$\chi_{{\phantom{}}_{\mu}} $ and $\chi_{{\phantom{}}_{\mu\nu}} $.

\head V. Appendix
\endhead

{\bf A.1.} In this part of the Appendix we give the explicit formulae for 
the  adjoint action of the $\kap$ on the elements $\vd^\a{}_\be$, 
$\vd^\mu{}_\nu$, $\vd^{\mu\nu\a}$ and $x^{\a\be}$
$$
\aligned
\adx (\vd^\a{}_\be \vd^\mu{}_\nu) &=  \vd^\rho{}_\si \vd^\gm{}_\de \otimes 
\vl_\rho{}^\mu \vl^\si {}_\nu \vl_\gm{}^\a \vl^\de{}_\be,\\
 \adx(\vd^{\mu\nu\a}) & =  \vd^\rho{}_\si  \vd^\be{}_\gm \otimes 
\vl_\rho{}^\mu \vl_\be{}^\a \vl^{\si\nu} x^\gm \\
& + \vd^{\rho\si\be} \otimes \vl_\rho{}^\mu \vl_\be{}^\a \vl_\si{}^\nu,\\
\adx(x^{\a\be}) & = x^{\mu\nu} \otimes \vl_\mu{}^\a \vl_\nu{}^\be + 
(\vd^{\mu\rho\nu} + \vd^{\nu\rho\mu}) \otimes  \vl_\mu{}^\a \vl_\nu{}^\be 
x_\rho\\
& + \vd^\mu{}_\rho \vd^\nu{}_\si \otimes  \vl_\mu{}^\a \vl_\nu{}^\be \Big( 
x^\rho x^\si - \jak g^{\rho\si} x^0 + \jak g^{0\si} x^\rho\Big).
\endaligned
\tag"{(A1)}"
$$
From (A1) we conclude that $(\vd^\a{}_\be \vd^\mu{}_\nu, 
\vd^{\mu\nu\a},x^{\mu\nu})$ span a linear ad-invariant set. Moreover, with 
respect to the Lorentz part of  \1 they transform as the corresponding 
tensors; note also that only the symmetric (with respect to $\mu$, $\nu$) 
part of $\vd^{\mu\a\nu}$ enters the transformation rule for $x^{\a\be}$. 
Therefore the linear set spanned by $\vd^\a{}_\be \vd ^\mu{}_\nu$,  
$\widetilde{\vd}^{\mu\nu\a}$ and $\widetilde{x}^{\mu\nu}$ which is obtained 
by  subtracting the completely antisymmetric subrepresentation from 
${\vd}^{\mu\nu\a}$ and the scalar (trace) subrepresentation from 
${x}^{\mu\nu}$ is also ad-invariant. By Lemma 1.7 of [12] the ideal $\car$ 
is ad-invariant.

{\bf A.2.} The action of $\si$ is
$$
\align
\si(\om^\mu{}_\nu \otimes \om^\a{}_\be) & = \om^\a{}_\be \otimes 
\om^\mu{}_\nu,\\
\si(\om^\mu{}_\nu \otimes \om) & = \om \otimes \om^\mu{}_\nu,\\
\si(\om^\mu \otimes \om) & = \om \otimes \om^\mu,
\endalign
$$
$$
\align
\si(\om \otimes \om) & = \om \otimes \om,\\
\si(\vo \otimes \om) & = \om \otimes \vo,\\
\si(\vo \otimes \vo) & = \vo \otimes \vo,\\
\si(\om \otimes \vo) & = \vo \otimes \om,\\
\si(\om^\mu{}_\nu \otimes \vo) & = \vo \otimes \om^\mu{}_\nu,\\
\si(\om^\mu \otimes \vo) & = \vo \otimes \om^\mu,\\
\si(\vo \otimes \om^\mu{}_\nu) & =  \om^\mu{}_\nu \otimes \vo ,\\
\si(\vo \otimes \om^\mu) & =  \om^\mu \otimes \vo - \dfrac{3}{\kappa^2} 
\ve_{\be\rho\nu}   \om^{\mu\be} \otimes \om^{\rho\nu} ,\\
\si(\om^\mu{}_\nu  \otimes \om^\a) & =  \om^\a \otimes \om^\mu{}_\nu + \jak 
(\de^0_\nu \om^\a{}_\rho \otimes \om^{\mu\rho} \\
& + \de^\mu{}_0 \om^\a{}_\rho \otimes \om^\rho{}_\nu - \om^\a{}_\nu \otimes 
\om^{\mu 0} -  \om^{\a\mu} \otimes  \om^0{}_\nu)\\
& + \dfrac{1}{6} \ve^\mu{}_\nu{}^\si \om^\a{}_\si \otimes \vo,\\ 
\si(\om^\a \otimes \om^\mu{}_\nu) & =  \om^\mu{}_\nu \otimes \om^\a 
+ \jak (\de^0_\nu \om^\mu{}_\si \otimes \om^{\si\a} + \om^\mu{}_0 \otimes 
\om^\a_{}\nu\\
& +  \de^\mu_0 \om^\si{}_\nu \otimes \om_\si{}^\a + \om_{0\nu} \otimes 
\om^{\a\mu}) + \dfrac{1}{6} \ve^\si{}_\nu{}^\a \om^\mu{}_\si \otimes \vo\\
& + \dfrac{1}{6} \ve^{\si\mu\a}  \om_{\si\nu}  \otimes \vo,\\
\si(\om \otimes \om^\mu) & =  \om^\mu  \otimes \om 
+ \dfrac{3i}{\kappa^3} ( \om^{\rho\mu} \otimes \om_{\rho 0} + \om_{\rho 0}
 \otimes \om^{\rho\mu} - \de^\mu_0 \om^\rho{}_\si  \otimes \om_\si{}^\rho)\\
& - \dfrac{3}{\kappa^2} ( \om^{\mu\rho} \otimes \om_{\rho} + \om_{\rho} \otimes 
\om^{\rho\mu} ) + \dfrac{1}{2\kappa^2} \ve^{\mu\rho\si} \om_{\rho\si} 
\otimes \vo,\\
\si(\om \otimes \om^\mu{}_\nu) & =  \om^\mu{}_\nu \otimes \om 
- \dfrac{3}{\kappa^2} ( \om^\si{}_\nu \otimes \om_\si{}^\mu - \om_\si{}^\mu 
 \otimes \om^{\si}{}_\nu)  ,\\
\si(\om^\mu \otimes \om^\nu) & = \om^\nu \otimes \om^\mu + \dfrac{1}{3} 
\om^{\nu\mu} \otimes \om + \dfrac{i}{\kappa} (\om^0 \otimes \om^{\mu\nu} +  
\om^{\mu\nu} \otimes \om^0)\\ 
& + \dfrac{1}{\kappa^2} ( \om_0{}^\nu \otimes \om^\mu{}_0 + \om_0{}^\mu 
\otimes  \om_0{}^\nu)\\
& + \dfrac{1}{\kappa^2} ( \de^\nu_0 \om^{\rho\mu} \otimes \om_{\rho 0} + 
\de^\mu_0 \om_{\rho 0} \otimes \om^{\rho\nu})\\
& + \jak  ( \de^\nu_0 \om^{\rho} \otimes \om_{\rho}{}^\mu + \de^\mu_0 
\om^\nu{}_\rho \otimes \om^{\rho}\\
& - \jakk \de^\nu_0 \de^\mu_0 \om^\rho{}_\si \otimes \om_\rho{}^\si - 
\dfrac{i}{6\kappa} \ve^{\nu\mu\si} \om_{0\si}  \otimes \vo\\ 
& - \dfrac{1}{6}  \ve^\nu{}_\rho{}^\mu \om^\rho  \otimes \vo \\
& + \dfrac{1}{2\kappa^2}  \ve^\nu{}_\rho{}^\mu \ve_{\si\tau\lm} 
\om^{\rho\si} \otimes \om^{\tau\lm}.
\endalign
$$
The action of $\si$ seems to be as much complicated as in four-\5al case 
[13]. In spite of that the exterior \3 appears to be simpler: the \5 of 
$\vg^{\wedge 2}$ equals $\binom{\dim \vg}{2}$.

\enddocument